

\input amssym.tex 

\def\unredoffs{}
\tolerance=1000\hfuzz=2pt
\catcode`\@=11 
\ifx\hyperdef\UNd@FiNeD\def\hyperdef#1#2#3#4{#4}\def\hyperref#1#2#3#4{#4}\def\href#1#2{#2}\fi
\magnification=1200\unredoffs\baselineskip=16pt plus 2pt minus 1pt
\def\Date#1{\vfill\leftline{#1}\tenpoint\supereject%
\footline={\hss\tenrm\hyperdef\hypernoname{page}\folio\folio\hss}}%

{\count255=\time\divide\count255 by 60 \xdef\hourmin{\number\count255}
 \multiply\count255 by-60\advance\count255 by\time
 \xdef\hourmin{\hourmin:\ifnum\count255<10 0\fi\the\count255}
}
\def\date{\number\day.\number\month.\number\year\ at \hourmin}


\def\nolabels{\def\wrlabeL##1{}\def\eqlabeL##1{}\def\reflabeL##1{}}
\def\writelabels{\def\wrlabeL##1{\leavevmode\vadjust{\rlap{\smash%
{\line{{\escapechar=` \hfill\rlap{\sevenrm\hskip.03in\string##1}}}}}}}%
\def\eqlabeL##1{{\escapechar-1\rlap{\sevenrm\hskip.05in\string##1}}}%
\def\reflabeL##1{\noexpand\llap{\noexpand\sevenrm\string\string\string##1}}}
\nolabels

\global\newcount\secno \global\secno=0
\global\newcount\meqno \global\meqno=1
\def\s@csym{}

\def\newsec#1\par{\global\advance\secno by1%
{\toks0{#1}\message{(\the\secno. \the\toks0)}}%
\global\subsecno=0\eqnres@t\let\s@csym\secsym\xdef\secn@m{\the\secno}\noindent
{\bf\hyperdef\hypernoname{section}{\the\secno}{\the\secno.} #1}%
\writetoca{{\string\hyperref{}{section}{\the\secno}{\bf \the\secno\quad}} {\bf #1}}\par%
\nobreak\medskip\nobreak\noindent\ignorespaces}
\def\eqnres@t{\xdef\secsym{\the\secno.}\global\meqno=1\bigbreak\bigskip}
\def\sequentialequations{\def\eqnres@t{\bigbreak}}\xdef\secsym{}

\global\newcount\subsecno \global\subsecno=0
\def\subsec#1\par{\global\advance\subsecno by1%
{\toks0{#1}\message{(\s@csym\the\subsecno. \the\toks0)}}%
\global\subsubsecno=0%
\ifnum\lastpenalty>9000\else\bigbreak\fi
\noindent{\it\hyperdef\hypernoname{subsection}{\secn@m.\the\subsecno}%
{\secn@m.\the\subsecno.} #1}\writetoca{\string\hskip1.45cm
{\string\hyperref{}{subsection}{\secn@m.\the\subsecno}{\secn@m.\the\subsecno.}}
{#1}}\par\nobreak\medskip\nobreak\noindent\ignorespaces}

\global\newcount\subsubsecno \global\subsubsecno=0
\def\subsubsec#1\par{\global\advance\subsubsecno by1%
{\toks0{#1}\message{(\secn@m.\the\subsecno.\the\subsubsecno. \the\toks0)}}%
\global\subsubsubsecno=0%
\ifnum\lastpenalty>9000\else\bigbreak\fi
\noindent{\it\hyperdef\hypernoname{subsubsection}{\secn@m.\the\subsecno\the\subsubsecno}%
{\secn@m.\the\subsecno.\the\subsubsecno.} #1}
\par\nobreak\medskip\nobreak\noindent\ignorespaces}

\global\newcount\subsubsubsecno \global\subsubsubsecno=0
\def\subsubsubsec#1\par{\global\advance\subsubsubsecno by1%
{\toks0{#1}\message{(\secn@m.\the\subsecno.\the\subsubsecno.\the\subsubsubsecno \the\toks0)}}%
\ifnum\lastpenalty>9000\else\bigbreak\fi
\noindent{\it\hyperdef\hypernoname{subsubsection}{\secn@m.\the\subsecno\the\subsubsecno\the\subsubsubsecno}%
{\secn@m.\the\subsecno.\the\subsubsecno.\the\subsubsubsecno.} #1}%
\par\nobreak\medskip\nobreak\noindent\ignorespaces}


\def\newnewsec#1#2\par{\global\advance\secno by1%
{\toks0{#2}\message{(\secn@m. \the\toks0)}}%
\global\subsecno=0\global\subsubsecno=0\eqnres@t\let\s@csym\secsym\xdef\secn@m{\the\secno}\noindent
\ifnum\lastpenalty>9000\else\bigbreak\fi
\noindent{\bf\hyperdef\hypernoname{section}{\secn@m}{\secn@m.} #2}%
\writetoca{{\string\hyperref{}{section}{\the\secno}{\bf \the\secno\quad}} {\bf #2}}
\DefWarn#1%
\xdef#1{\noexpand\hyperref{}{section}{\the\secno}%
{\the\secno}}\writedef{#1\leftbracket#1}\wrlabeL{#1=#1}%
\par\nobreak\medskip\nobreak\noindent\ignorespaces}

\def\newsubsec#1#2\par{\global\advance\subsecno by1%
{\toks0{#2}\message{(\secn@m.\the\subsecno. \the\toks0)}}%
\global\subsubsecno=0%
\ifnum\lastpenalty>9000\else\bigbreak\fi
\noindent{\it\hyperdef\hypernoname{subsection}{\secn@m.\the\subsecno}%
{\secn@m.\the\subsecno.} #2}
\DefWarn#1%
\xdef#1{\noexpand\hyperref{}{subsection}{\secn@m.\the\subsecno}%
{\secn@m.\the\subsecno}}\writedef{#1\leftbracket#1}\wrlabeL{#1=#1}%
\writetoca{\string\hskip1.45cm
{\string\hyperref{}{subsection}{\secn@m.\the\subsecno}{\secn@m.\the\subsecno.}}
{#2}}%
\par\nobreak\medskip\nobreak\noindent\ignorespaces}

\def\newsubsubsec#1#2\par{\global\advance\subsubsecno by1%
{\toks0{#2}\message{(\secn@m.\the\subsecno.\the\subsubsecno. \the\toks0)}}%
\global\subsubsubsecno=0%
\ifnum\lastpenalty>9000\else\bigbreak\fi
\noindent{\it\hyperdef\hypernoname{subsubsection}{\secn@m.\the\subsecno\the\subsubsecno}%
{\secn@m.\the\subsecno.\the\subsubsecno.} #2}
\DefWarn#1%
\xdef#1{\noexpand\hyperref{}{subsubsection}{\secn@m.\the\subsecno.\the\subsubsecno}%
{\secn@m.\the\subsecno.\the\subsubsecno}}\writedef{#1\leftbracket#1}\wrlabeL{#1=#1}%
\par\nobreak\medskip\nobreak\noindent\ignorespaces}

\def\newsubsubsubsec#1#2\par{\global\advance\subsubsubsecno by1%
{\toks0{#2}\message{(\secn@m.\the\subsecno.\the\subsubsecno.\the\subsubsubsecno \the\toks0)}}%
\ifnum\lastpenalty>9000\else\bigbreak\fi
\noindent{\it\hyperdef\hypernoname{subsubsection}{\secn@m.\the\subsecno\the\subsubsecno\the\subsubsubsecno}%
{\secn@m.\the\subsecno.\the\subsubsecno.\the\subsubsubsecno.} #2}
\DefWarn#1%
\xdef#1{\noexpand\hyperref{}{subsubsubsection}{\secn@m.\the\subsecno.\the\subsubsecno.\the\subsubsubsecno}%
{\secn@m.\the\subsecno.\the\subsubsecno.\the\subsubsubsecno}}\writedef{#1\leftbracket#1}\wrlabeL{#1=#1}%
\par\nobreak\medskip\nobreak\noindent\ignorespaces}

\def\appendix#1#2{\global\meqno=1\global\subsecno=0\global\subsubsecno=0\xdef\secsym{\hbox{#1.}}%
\bigbreak\bigskip\noindent{\bf Appendix \hyperdef\hypernoname{appendix}{#1}%
{#1.} #2}{\toks0{(#1. #2)}\message{\the\toks0}}%
\xdef\s@csym{#1.}\xdef\secn@m{#1}%
\writetoca{{\string\hyperref{}{appendix}{#1}{\bf {#1}\quad}} {\bf #2}}%
\par\nobreak\medskip\nobreak}

%
\def\checkm@de#1#2{\ifmmode{\def\f@rst##1{##1}\hyperdef\hypernoname{equation}%
{#1}{#2}}\else\hyperref{}{equation}{#1}{#2}\fi}
\def\eqnn#1{\DefWarn#1\xdef #1{(\noexpand\relax\noexpand\checkm@de%
{\s@csym\the\meqno}{\secsym\the\meqno})}%
\wrlabeL#1\writedef{#1\leftbracket#1}\global\advance\meqno by1}
\def\f@rst#1{\c@t#1a\em@ark}\def\c@t#1#2\em@ark{#1}
\def\eqna#1{\DefWarn#1\wrlabeL{#1$\{\}$}%
\xdef #1##1{(\noexpand\relax\noexpand\checkm@de%
{\s@csym\the\meqno\noexpand\f@rst{##1}1}{\hbox{$\secsym\the\meqno##1$}})}
\writedef{#1\numbersign1\leftbracket#1{\numbersign1}}\global\advance\meqno by1}
\def\eqn#1#2{\DefWarn#1%
\xdef #1{(\noexpand\hyperref{}{equation}{\s@csym\the\meqno}%
{\secsym\the\meqno})}$$#2\eqno(\hyperdef\hypernoname{equation}%
{\s@csym\the\meqno}{\secsym\the\meqno})\eqlabeL#1$$%
\writedef{#1\leftbracket#1}\global\advance\meqno by1}
\def\xeqn{\expandafter\xe@n}\def\xe@n(#1){#1}
\def\xeqna#1{\expandafter\xe@n#1}
\def\eqns#1{(\e@ns #1{\hbox{}})}
\def\e@ns#1{\ifx\UNd@FiNeD#1\message{eqnlabel \string#1 is undefined.}%
\xdef#1{(?.?)}\fi{\let\hyperref=\relax\xdef\next{#1}}%
\ifx\next\em@rk\def\next{}\else%
\ifx\next#1\xeqn#1\else\def\n@xt{#1}\ifx\n@xt\next#1\else\xeqna#1\fi
\fi\let\next=\e@ns\fi\next}
\def\DefWarn#1{}
%
\newskip\footskip\footskip14pt plus 1pt minus 1pt 
\def\footnotefont{\ninepoint}\def\f@t#1{\footnotefont #1\@foot}
\def\f@@t{\baselineskip\footskip\bgroup\footnotefont\aftergroup\@foot\let\next}
\setbox\strutbox=\hbox{\vrule height9.5pt depth4.5pt width0pt}
\global\newcount\ftno \global\ftno=0
\def\foot{\global\advance\ftno by1\def\foot@rg{\hyperref{}{footnote}%
{\the\ftno}{\the\ftno}\xdef\foot@rg{\noexpand\hyperdef\noexpand\hypernoname%
{footnote}{\the\ftno}{\the\ftno}}}\footnote{$^{\foot@rg}$}}
%
%
%
\global\newcount\refno \global\refno=1
\newwrite\rfile
\def\ref{[\hyperref{}{reference}{\the\refno}{\the\refno}]\nref}
\def\nref#1{\DefWarn#1%
\xdef#1{[\noexpand\hyperref{}{reference}{\the\refno}{\the\refno}]}%
\writedef{#1\leftbracket#1}%
\ifnum\refno=1\immediate\openout\rfile=\jobname.refs\fi
\chardef\wfile=\rfile\immediate\write\rfile{\noexpand\item{[\noexpand\hyperdef%
\noexpand\hypernoname{reference}{\the\refno}{\the\refno}]\ }%
\reflabeL{#1\hskip.31in}\pctsign}\global\advance\refno by1\findarg}
\def\findarg#1#{\begingroup\obeylines\newlinechar=`\^^M\pass@rg}
{\obeylines\gdef\pass@rg#1{\writ@line\relax #1^^M\hbox{}^^M}%
\gdef\writ@line#1^^M{\expandafter\toks0\expandafter{\striprel@x #1}%
\edef\next{\the\toks0}\ifx\next\em@rk\let\next=\endgroup\else\ifx\next\empty%
\else\immediate\write\wfile{\the\toks0}\fi\let\next=\writ@line\fi\next\relax}}
\def\striprel@x#1{} \def\em@rk{\hbox{}}
\def\lref{\begingroup\obeylines\lr@f}
\def\lr@f#1#2{\DefWarn#1\gdef#1{\let#1=\UNd@FiNeD\ref#1{#2}}\endgroup\unskip}
\def\semi{;\hfil\break}
\def\addref#1{\immediate\write\rfile{\noexpand\item{}#1}} 
\def\listrefs{\vfill\supereject\immediate\closeout\rfile\writestoppt
\baselineskip=\footskip\centerline{{\bf References}}\bigskip{\parindent=20pt%
\frenchspacing\escapechar=` \input \jobname.refs\vfill\eject}\nonfrenchspacing}
\def\startrefs#1{\immediate\openout\rfile=\jobname.refs\refno=#1}
\def\xref{\expandafter\xr@f}\def\xr@f[#1]{#1}
\def\refs#1{\count255=1[\r@fs #1{\hbox{}}]}
\def\r@fs#1{\ifx\UNd@FiNeD#1\message{reflabel \string#1 is undefined.}%
\nref#1{need to supply reference \string#1.}\fi%
\vphantom{\hphantom{#1}}{\let\hyperref=\relax\xdef\next{#1}}%
\ifx\next\em@rk\def\next{}%
\else\ifx\next#1\ifodd\count255\relax\xref#1\count255=0\fi%
\else#1\count255=1\fi\let\next=\r@fs\fi\next}
%

%
\newwrite\ffile\global\newcount\figno \global\figno=1
\def\fig{fig.~\hyperref{}{figure}{\the\figno}{\the\figno}\nfig}
\def\nfig#1{\DefWarn#1%
\xdef#1{fig.~\noexpand\hyperref{}{figure}{\the\figno}{\the\figno}}%
\writedef{#1\leftbracket fig.\noexpand~\xfig#1}%
\ifnum\figno=1\immediate\openout\ffile=\jobname.figs\fi\chardef\wfile=\ffile%
{\let\hyperref=\relax
\immediate\write\ffile{\noexpand\medskip\noexpand\item{Fig.\ %
\noexpand\hyperdef\noexpand\hypernoname{figure}{\the\figno}{\the\figno}. }
\reflabeL{#1\hskip.55in}\pctsign}}\global\advance\figno by1\findarg}
\def\xfig{\expandafter\xf@g}\def\xf@g fig.\penalty\@M\ {}
\def\figs#1{figs.~\f@gs #1{\hbox{}}}
\def\f@gs#1{{\let\hyperref=\relax\xdef\next{#1}}\ifx\next\em@rk\def\next{}\else
\ifx\next#1\xfig #1\else#1\fi\let\next=\f@gs\fi\next}
%
\def\figin{\epsfcheck\figin}\def\figins{\epsfcheck\figins}
\def\epsfcheck{\ifx\epsfbox\UnDeFiNeD
\message{(NO epsf.tex, FIGURES WILL BE IGNORED)}
\gdef\figin##1{\vskip2in}\gdef\figins##1{\hskip.5in}
\else\message{(FIGURES WILL BE INCLUDED)}%
\gdef\figin##1{##1}\gdef\figins##1{##1}\fi}
\def\figinsert{\goodbreak\topinsert}
\def\ifig#1#2#3{\DefWarn#1\xdef#1{fig.~\the\figno}
\writedef{#1\leftbracket fig.\noexpand~\the\figno}%
\figinsert\figin{\centerline{#3}}
\smallskip
\leftskip=0pt \rightskip=0pt
\baselineskip12pt\noindent
{{\bf Fig.~\the\figno}\ \ninepoint #2}
\medskip
\global\advance\figno by1\par\endinsert}
\newwrite\lfile
{\escapechar-1\xdef\pctsign{\string\%}\xdef\leftbracket{\string\{}
\xdef\rightbracket{\string\}}\xdef\numbersign{\string\#}}
\def\writedefs{\immediate\openout\lfile=label.defs \def\writedef##1{%
{\let\hyperref=\relax\let\hyperdef=\relax\let\hypernoname=\relax
 \immediate\write\lfile{\string\checkdef\string##1\rightbracket}}}}%
\def\writestop{\def\writestoppt{\immediate\write\lfile{\string\pageno
 \the\pageno\string\startrefs\leftbracket\the\refno\rightbracket
 \string\def\string\secsym\leftbracket\secsym\rightbracket
 \string\secno\the\secno\string\meqno\the\meqno}\immediate\closeout\lfile}}
\def\writestoppt{}\def\writedef#1{}

\def\seclab#1\par{\DefWarn#1%
\xdef #1{\noexpand\hyperref{}{section}{\the\secno}{\the\secno}}%
\writedef{#1\leftbracket#1}\wrlabeL{#1=#1}\par%
\nobreak\medskip\nobreak\noindent\ignorespaces}
\def\subseclab#1\par{\DefWarn#1%
\xdef #1{\noexpand\hyperref{}{subsection}{\the\secno.\the\subsecno}%
{\the\secno.\the\subsecno}}\writedef{#1\leftbracket#1}\wrlabeL{#1=#1}\par%
\nobreak\medskip\nobreak\noindent\ignorespaces}
\def\subsubseclab#1\par{\DefWarn#1%
\xdef#1{\noexpand\hyperref{}{subsubsection}{\the\secno.\the\subsecno.\the\subsubsecno}%
{\the\secno.\the\subsecno.\the\subsubsecno}}\writedef{#1\leftbracket#1}\wrlabeL{#1=#1}\par%
\nobreak\medskip\nobreak\noindent\ignorespaces}
\def\applab#1\par{\DefWarn#1%
\xdef#1{\noexpand\hyperref{}{appendix}{\secn@m}{\secn@m}}%
\writedef{#1\leftbracket#1}\wrlabeL{#1=#1}%
\par\nobreak\medskip\nobreak\noindent\ignorespaces}
\def\appsublab#1{\DefWarn#1%
\xdef #1{\noexpand\hyperref{}{appendix}{\secn@m.\the\subsecno}{\secn@m.\the\subsecno}}%
\writedef{#1\leftbracket#1}\wrlabeL{#1=#1}}
\newwrite\tfile \def\writetoca#1{}
\def\leaderfill{\leaders\hbox to 1em{\hss.\hss}\hfill}
\def\writetoc{\immediate\openout\tfile=\jobname.toc
   \def\writetoca##1{{\edef\next{\write\tfile{\noindent ##1
   \string\leaderfill{
   \string\hyperref{}{page}{\noexpand\number\pageno}%
   {\noexpand\number\pageno}} \par}}\next}}
}
\newread\ch@ckfile
\def\listtoc{\immediate\closeout\tfile\immediate\openin\ch@ckfile=\jobname.toc
\ifeof\ch@ckfile\message{no file \jobname.toc, no table of contents this pass}%
\else\closein\ch@ckfile\centerline{\bf Contents}\nobreak\medskip%
{\baselineskip=15.5pt\footnotefont\parskip=0pt\catcode`\@=11\input\jobname.toc
\catcode`\@=12\bigbreak\bigskip}\fi}
\catcode`\@=12 
\def\tenpoint{\def\rm{\fam0\tenrm}
\textfont0=\tenrm \scriptfont0=\sevenrm \scriptscriptfont0=\fiverm
\textfont1=\teni  \scriptfont1=\seveni  \scriptscriptfont1=\fivei
\textfont2=\tensy \scriptfont2=\sevensy \scriptscriptfont2=\fivesy
\textfont\itfam=\tenit \def\it{\fam\itfam\tenit}\def\footnotefont{\ninepoint}%
\textfont\bffam=\tenbf \def\bf{\fam\bffam\tenbf}\def\sl{\fam\slfam\tensl}\rm}
\font\ninerm=cmr9 \font\sixrm=cmr6 \font\ninei=cmmi9 \font\sixi=cmmi6
\font\ninesy=cmsy9 \font\sixsy=cmsy6 \font\ninebf=cmbx9
\font\nineit=cmti9 \font\ninesl=cmsl9 \skewchar\ninei='177
\skewchar\sixi='177 \skewchar\ninesy='60 \skewchar\sixsy='60
\def\ninepoint{\def\rm{\fam0\ninerm}
\textfont0=\ninerm \scriptfont0=\sixrm \scriptscriptfont0=\fiverm
\textfont1=\ninei \scriptfont1=\sixi \scriptscriptfont1=\fivei
\textfont2=\ninesy \scriptfont2=\sixsy \scriptscriptfont2=\fivesy
\textfont\itfam=\ninei \def\it{\fam\itfam\nineit}\def\sl{\fam\slfam\ninesl}%
\textfont\bffam=\ninebf \def\bf{\fam\bffam\ninebf}\rm}
%
\hyphenation{anom-aly anom-alies coun-ter-term coun-ter-terms}

\def\tikzcaption#1#2{\DefWarn#1\xdef#1{Fig.~\the\figno}
\writedef{#1\leftbracket Fig.\noexpand~\the\figno}%
{
\smallskip
\leftskip=20pt \rightskip=20pt \baselineskip12pt\noindent
{{\bf Fig.~\the\figno}\ \ninepoint #2}
\bigskip
\global\advance\figno by1 \par}}

\def\ntoalpha#1{%
\ifcase#1%
@%
\or A\or B\or C\or D\or E\or F\or G\or H\or I\or J\or K\or L\or M%
\fi
}

\global\newcount\appno \global\appno=1
\def\applab#1{\xdef #1{\ntoalpha{\appno}}\writedef{#1\leftbracket#1}\wrlabeL{#1=#1}
\global\advance\appno by1}

\def\preprint#1 #2\par{\rightline{\vbox{\baselineskip12pt\hbox{#1}\hbox{#2}}}\vskip2cm}
%
\def\title#1\par{\centerline{\bf #1}\nopagenumbers\pageno=0}
\def\author#1\par{\bigskip\bigskip\centerline{#1}}

\newcount\addressno

\def\email#1#2{
\footnote{\null}{\kern-\parindent \llap{$^#1$\hskip1pt}email: #2}}

\def\startcenter{%
  \par
  \begingroup
  \leftskip=0pt plus 1fil
  \rightskip=\leftskip
  \parindent=0pt
  \parfillskip=0pt
}
\def\stopcenter{\endgroup}

\def\address{\bigskip%
  \ifnum\the\addressno=0\else\stopcenter\endgroup\fi
  \advance\addressno by 1%
  \begingroup
  \startcenter
  \it
  \obeylines
  \addressAux
}
\def\addressAux#1{#1}

\def\abstract{\stopcenter\endgroup\bigskip\bigskip\noindent}

\def\Dsl{\,\raise.15ex\hbox{/}\mkern-13.5mu D} 
\def\dsl{\raise.15ex\hbox{/}\kern-.57em\partial}
 \def\Tr{{\rm Tr}}
\def\boxeqn#1{\vcenter{\vbox{\hrule\hbox{\vrule\kern3pt\vbox{\kern3pt
	\hbox{${\displaystyle #1}$}\kern3pt}\kern3pt\vrule}\hrule}}}


\def\d{{\delta}}

\def\s{{\sigma}}

\def\half{{1\over 2}}

\def\bar{\overline}
\def\({\left(}
\def\){\right)}



\def\len#1{{%
\def\Dlen{\left|\mkern-1mu #1\mkern -0.5mu\right|}%
\def\Sslen{\left|\mkern-1.3mu #1\mkern -1.3mu\right|}%
\def\SSlen{\left|\mkern-2.8mu #1\mkern-1.3mu\right|}%
\mathchoice{\Dlen}{\Dlen}{\Sslen}{\SSlen}}}

\def\sfrac#1/#2{\kern.1em\raise.5ex\hbox{\the\scriptfont0 #1}%
\kern-.1em/\kern-.15em\lower.25ex\hbox{\the\scriptfont0 #2}}

\font\tenshuffle=shuffle10 \font\sevenshuffle=shuffle7 \font\fiveshuffle=shuffle7 at 5pt
\def\shuffle{{%
\def\Dshuffle{\mathbin{\hbox{\tenshuffle\char'001}}}%
\def\Sshuffle{\mathbin{\hbox{\sevenshuffle\char'001}}}%
\def\SSshuffle{\mathbin{\hbox{\fiveshuffle\char'001}}}%
\mathchoice{\Dshuffle}{\Dshuffle}{\Sshuffle}{\SSshuffle}}}


\def\qed{\hbox{\hskip 3pt
\vbox{\hrule\hbox to 7pt{\vrule height 7pt\hfill\vrule}
\hrule}}\hskip3pt}

\overfullrule=0pt\relax

\frenchspacing

\def\checkdef#1#2{
\ifx\UndeFined#1%
	\def#1{#2}
\else
	\immediate\write16{*** BUG ***: the label \string#1 is already defined ***}
\fi
}
\newread\instream

\openin\instream= label.defs
\ifeof\instream\message{No labels in advance yet. Wait till next pass.}
\else\closein\instream \input label.defs
\fi
\writedefs

\def\arXiv:#1].{\hepthStrip#1 \nil}
\def\hepthStrip#1 #2\nil{\href{http://arxiv.org/abs/#1}{arXiv:#1 #2\unskip}].}


\def\Tr{{\rm Tr}}
\def\frac#1#2{{#1 \over #2}}

\title A closed-formula solution to the color-trace decomposition problem

\author
Ruggero Bandiera\email{\star}{bandiera@mat.uniroma1.it}$^{\star}$ and
Carlos R. Mafra\email{\dagger}{c.r.mafra@soton.ac.uk}$^{\dagger}$

\address
$^\star$ Universit\`a degli studi di Roma La Sapienza,
Dipartimento di Matematica ``Guido Castelnuovo",
P.le Aldo Moro 5, I-00185 Roma, Italy.

\address
$^\dagger$Mathematical Sciences and STAG Research Centre, University of Southampton,
Highfield, Southampton, SO17 1BJ, UK

\abstract
In these notes we present a closed-formula solution to the problem of decomposing traces of Lie algebra
generators into symmetrized traces and structure constants. The solution is written in terms of
Solomon idempotents and exploits a projection derived by Solomon in his work on the Poincar\'e-Birkhoff-Witt
theorem.

\Date{September 2020}


\lref\oneloopbb{
	C.R.~Mafra and O.~Schlotterer,
  	``The Structure of n-Point One-Loop Open Superstring Amplitudes,''
	JHEP {\bf 1408}, 099 (2014).
	[arXiv:1203.6215 [hep-th]].
}
\lref\Vcolor{
	T.~van Ritbergen, A.~N.~Schellekens and J.~A.~M.~Vermaseren,
  	``Group theory factors for Feynman diagrams,''
	Int.\ J.\ Mod.\ Phys.\ A {\bf 14}, 41 (1999).
	[hep-ph/9802376].
}
\lref\ddm{
	V.~Del Duca, L.~J.~Dixon and F.~Maltoni,
  	``New color decompositions for gauge amplitudes at tree and loop level,''
	Nucl.\ Phys.\ B {\bf 571}, 51 (2000).
	[hep-ph/9910563].
}

\lref\Reutenauer{
	C.~Reutenauer,
	``Free Lie Algebras,''
	London Mathematical Society Monographs, 1993
}
\lref\PBWReutenauer{
	C.~Reutenauer,
	``Theorem of Poincar\'e-Birkhoff-Witt, logarithm and
	symmetric group representations of degrees equal to
	Stirling numbers''. In Combinatoire \'enum\'erative (pp. 267-284) 1986. Springer, Berlin, Heidelberg.
}
\lref\Ree{
	R. Ree, ``Lie elements and an algebra associated with shuffles'',
	Ann. Math. {\bf 62}, No. 2 (1958), 210--220.
}
\lref\bandiera{
	R.~Bandiera, F.~Schaetz, ``Eulerian idempotent, pre-Lie logarithm and combinatorics of trees''.
	arXiv:1702.08907.
}
\lref\Loday{
	J.L.~Loday, ``S\'erie de Hausdorff, idempotents Eul\'eriens et algebres de Hopf''.
	Exp. Math. {\bf 12} (1994), 165-178.
}
\lref\solomon{
	L.~Solomon ``On the Poincar\'e-Birkhoff-Witt theorem''.
	Journal of Combinatorial Theory. 1968 May 1;4(4):363-75.
}
\lref\FORM{
	J.A.M.~Vermaseren,
	``New features of FORM,''
	arXiv:math-ph/0010025.
\semi
	M.~Tentyukov and J.A.M.~Vermaseren,
	``The multithreaded version of FORM,''
	arXiv:hep-ph/0702279.
}
\lref\gersten{
	M.~Gerstenhaber, ``Developments from Barr's thesis.''
	Journal of Pure and Applied Algebra 143, no. 1-3 (1999): 205-220.
}
\lref\giaquinto{
	A.~Giaquinto, ``Topics in algebraic deformation theory.''
	In Higher structures in geometry and physics, pp. 1-24. Birkh\"auser, Boston, MA, 2011.
}
\lref\garsia{
	A.M. Garsia, ``Combinatorics of the Free Lie Algebra and the Symmetric Group'',
	In Analysis, et Cetera, edited by Paul H. Rabinowitz and Eduard Zehnder,
	Academic Press, (1990) 309-382
}
\lref\casas{
	Arnal, Ana, Fernando Casas, and Cristina Chiralt.
	``A general formula for the Magnus expansion in terms of iterated integrals of right-nested commutators.''
	Journal of Physics Communications 2.3 (2018): 035024.
}



\newsec Introduction

The purpose of these notes is to present a closed-formula solution
to one of the problems addressed in \Vcolor\ via computer algebra \FORM. Given a simple Lie algebra whose generators $T^{a}$ satisfy
\eqn\commut{
[T^a,T^b] = i f^{abc} T^c\,,\quad \Tr(T^a T^b) = \half\delta^{ab}\,,
}
where $f^{abc}$ denote the totally anti-symmetric structure constants and
$\d^{ab}$ is the Kronecker delta, the problem
consists in expressing traces $\Tr(T^{a_1}\cdots T^{a_n})$ of products of Lie algebra generators $T^a$ (color factors) in terms of symmetrized traces 
\eqn\introsymtrace{ d^{a_1\cdots a_n} := \frac{1}{n!} \sum_{\sigma\in S_n}\Tr(T^{a_{\sigma(1)}}\cdots T^{a_{\sigma(n)}}) }
and structure constants $f^{abc}$. For example,
\eqn\introex{
\Tr(T^{a_1}T^{a_2}T^{a_3}) = d^{a_1a_2a_3} + {i\over 4}f^{a_1a_2a_3}\,.
}
These decompositions have important applications in the evaluation of loop amplitudes in perturbative field and string theories, as
they allow an efficient handling of their associated color structures in a manner described in \Vcolor.

While in \Vcolor\ an algorithm was obtained to generate these decompositions using computer algebra,
we will see here that this {\it color trace decomposition} problem admits an elegant closed-formula solution using
a result known in the free Lie algebra literature. The formula involves the so-called
{\it Solomon idempotent} or {\it first Eulerian idempotent}
\refs{\solomon,\garsia,\PBWReutenauer,\Loday,\bandiera} and its first few cases
are given by (to avoid cluttering we write $j$ instead of $a_j$)
\eqnn\elegant
$$\eqalignno{
\Tr(T^1 T^2) &= d^{12}\,, &\elegant\cr
\Tr(T^1 T^2 T^3) &= d^{123} +  d^{1a} E^{23}_a\,,\cr
\Tr(T^1 T^2 T^3 T^4) &=d^{1234}   +  d^{12a} E^{34}_a + d^{13a} E^{24}_a
+d^{14a} E^{23}_a +  d^{1a}E^{234}_a\,,\cr
\Tr(T^1 T^2 T^3 T^4 T^5) &=
d^{12345}
+ d^{123a} E^{45}_a
+ d^{124a}E^{35}_a
+ d^{125a} E^{34}_a
+d^{134a} E^{25}_a
+d^{135a} E^{24}_a\cr
&+d^{145a} E^{23}_a
+ d^{12a} E^{345}_a
+ d^{13a} E^{245}_a
+d^{14a} E^{235}_a
+d^{15a} E^{234}_a\cr
&+ d^{1ab} (E^{23}_a E^{45}_b
+ E^{24}_a E^{35}_b
+ E^{25}_a E^{34}_b)
+ d^{1a} E^{2345}_a\,,
}$$
where $E^{12\ldots n}_a$ denote the expansion coefficients of the Solomon idempotent with respect to the Lie algebra generators
$E(T^1\ldots T^n)=E^{1\ldots n}_aT^a$.  We refer to subsection  \eid\ for more precise definitions, and for now just point out
that these coefficients can be explicitly computed as polynomials in the structure constants $f^{abc}$, using results from
\bandiera\ (see also \casas). For instance, this yields the following solution to the color trace decomposition problem up to
$n=5$, 
\eqnn\Vtracetwo
$$\eqalignno{
\Tr(T^1 T^2) & = d^{12} = \half \d^{12}&\Vtracetwo\cr
\Tr(T^1T^2T^3) &= d^{123} + {i\over 4}f^{123}\cr
\Tr(T^1T^2T^3T^4) &=
       d^{1234} -{1\over6}f^{23a}f^{a41}+{1\over12}f^{24a}f^{a31}\cr
       &
       + {i\over2} d^{12a}f^{a34}
       + {i\over2} d^{13a}f^{a24}
       + {i\over2} d^{14a}f^{a23}\cr
\Tr(T^1 T^2T^3T^4 T^5) &=
        d^{12345}
         + {i\over24}\big(-3 f^{23 a}f^{a 4 b}f^{b 51}
         + f^{23 a}f^{a 5 b}f^{b 41}
         + f^{24 a}f^{a 3 b}f^{b 51} \cr
        	&\hskip50pt{}
         + f^{24 a}f^{a 5 b}f^{b 31}
         + f^{25 a}f^{a 3 b}f^{b 41}
         - f^{25 a}f^{a 4 b}f^{b 31}\big)\cr
        	&
       	 - {1\over4} d^{1ab} \big(
       	 f^{23 a} f^{45 b}
        + f^{24 a} f^{35 b}
        + f^{25 a} f^{34 b}
       	\big)\cr
       	&+ d^{12a} (-{1\over3}f^{34b}f^{b5a}+{1\over6}f^{35b}f^{b4a})
       	 + d^{13a} (-{1\over3}f^{24b}f^{b5a}+{1\over6}f^{25b}f^{b4a}) \cr
       	&+ d^{14a} (-{1\over3}f^{23b}f^{b5a}+{1\over6}f^{25b}f^{b3a})
       	 + d^{15a} (-{1\over3}f^{23b}f^{b4a}+{1\over6}f^{24b}f^{b3a})\cr
       	&
        + {i\over2}\big(
	 d^{123 a}f^{a 45}
       	+ d^{124 a}f^{a 35}
	+ d^{125 a}f^{a 34}
	+ d^{134 a}f^{a 25}
	+ d^{135 a}f^{a 24}
	+ d^{145 a}f^{a 23}
       	\big)\,,
}$$
recovering computations from \oneloopbb\foot{
The expansions in \oneloopbb\ use a different basis of color factors. That particular basis follows from 
expanding the Eulerian idempotents $E(x_1, \ldots,x_n)$ \eulerianDynkin\ in terms of the right-to-left free Lie algebra basis
$r(\ldots,x_n)$ rather than the left-to-right
$\ell(x_1, \ldots)$ as chosen in this work.}.

Our solution to the color trace decomposition problem shall depend on a projection formula due to Solomon \solomon, and related 
to the Poincar\'e-Birkhoff-Witt Theorem (see \Reutenauer). Recall that in particular the latter implies that a product of
generators $T^{p_1}\cdots T^{p_n}$ can be expanded as a linear combination of symmetrized products of Lie monomials in the
generators $T^{p_1},\ldots\,, T^{p_n}$ (for instance $T^1T^2=\frac{1}{2}(T^1T^2+T^2T^1) + \frac{1}{2}[T^1,T^2]$). Solomon's
formula provides such an expansion explicitly in terms of the first Eulerian idempotent, and from this (and the usual cyclic
properties of the trace) we shall deduce the following compact formula \closedshuffle, containing \elegant\ as particular cases.
Given a word $P=p_1\cdots p_n$, we denote by $T^P:=T^{p_1}\cdots T^{p_n}$ and by $\bar{\d}_k(P)=\sum_{(P)}
P_{(1)}\otimes\cdots\otimes P_{(k)}$ the $k$-th (reduced) deshuffle map applied to $P$, using Sweedler's notation (see the
subsection below for more precise definitions): then our formula reads
\eqn\closedshuffle{
\Tr(T^1T^P)
= \sum_{k\ge 1}\sum_{(P)}
{1\over k!}d^{1a_1a_2 \ldots a_k}E_{a_1}^{P_{(1)}}E_{a_2}^{P_{(2)}} \ldots E_{a_k}^{P_{(k)}}\,,
}
(where the second summation runs over the set of $k$-deshuffles $P_{(1)}\otimes\cdots\otimes P_{(k)}$ of $P$).  After expanding
the $E^{i_1\cdots i_k}_a$ as polynomials in the structure constants, we finally obtain the following closed formula solution for
the color trace decomposition problem
\eqn\eqclosedformula{
\Tr(T^0T^1\cdots T^n) = \sum_{S_n\ni \sigma = \sigma_1\cdots
\sigma_k} i^{n-k} C_{\sigma_1} \cdots C_{\sigma_k} d^{0a_1\cdots a_k} F^{\sigma_1}_{a_1}\cdots F^{\sigma_k}_{a_k}.
}
The latter formula deserves some explanations. First of all, the sum runs over the set of permutations $\sigma\in S_n$, which
are identified with the corresponding words $\sigma(1)\cdots\sigma(n)$. Then $\sigma=\sigma_1\cdots\sigma_k$ denotes the {\it
standard factorization} of $\sigma$, i.e., the unique factorization of $\sigma$ as the concatenation product of subwords
$\sigma_1,\ldots ,\sigma_k$ such that $\sigma_1>\cdots>\sigma_k$ in the lexicographical order and for all $1\leq j\leq k$ the
first letter in $\sigma_j$ is the minimum among its letters. A few examples are given by,
\eqn\exastdfac{
1432 = (1432)\,,\quad 2134 = (2)(134)\,,\quad 54132 = (5)(4)(132)\,,\quad 42671835 = (4)(267)(1835)\,.
}
Finally, given a word $P= p_1\cdots p_i$ we denote by $C_P:=
\frac{(-1)^{d_P}}{\len{P}{\len{P}-1 \choose d_P}}$ (where $d_P$ is the number of descents in $P$, once again 
we refer to the subsection below for more details) and by $F^{1}_a := \delta^{1a}$, $F^{12}_a
:= f^{12a}$, $F^{123}_a := f^{12b}f^{b3a}$, $F^{1234}_a:= f^{12c}f^{c3b}f^{b4a}$, in general,
\eqn\capF{
F^P_a=f^{p_1\cdots p_i}_a := f^{p_1p_2c_1}f^{c_1p_3c_2}\cdots f^{c_{i-2}p_ia}.
}
To better
understand the above formula \eqclosedformula\ the reader might check that for $n\le 4$ it precisely recovers \Vtracetwo\ (after
the obvious shift of indices). For instance, for $n=4$, $\sigma=2413\in S_4$, the standard factorization is
$\sigma=\sigma_1\sigma_2=(24)(13)$ and the corresponding term in \eqclosedformula \ is
$i^2C_{24}C_{13}d^{0ab}F^{24}_aF^{13}_a=-\frac{1}{4}d^{0ab}f^{24a}f^{13b}$. As further examples, for
$\sigma=3142=\sigma_1\sigma_2=(3)(142)$ we get the term $i^2C_3C_{142}d^{0ab}F^3_aF^{142}_b= \frac{1}{6}d^{03b}f^{14c}f^{c2b}$,
for $\sigma=4231=\sigma_1\sigma_2\sigma_3=(4)(23)(1)$ the one $iC_4C_{23}C_1d^{0abc}F^4_aF^{23}_bF^1_c =
\frac{i}{2}d^{014b}f^{b23}$ and for $\sigma = 1432 =\sigma_1$ the one $i^3C_{1432}d^{0a}F^{1432}_a = -\frac{i}{12}
d^{0a}f^{14c}f^{c3b}f^{b2a}= -\frac{i}{24} f^{02b}f^{b3c}f^{c14}$.  It is important to observe that the output of
\eqclosedformula\ is already written down in a basis of color factors, that is, no linear relations among its terms can be
deduced using only the Jacobi identities \genJacobi.

\newsubsec\notsec Notation on words

In this paper the labels in indices such as $a_2$ will be interpreted as {\it letters} from the alphabet of natural numbers
$\{1,2,3,\ldots \}$ and denoted by lower-case letters (e.g. $j=2$). {\it Words} composed of such letters will be denoted by
capital letters such as $P=13245$. The length of the word $P$ is denoted $|P|$ and it is given by the number of its letters.
Given a word $P$, a descent in $P$ is a pair of consecutive letters $P = \cdots p_jp_{j+1}\cdots$ such that $p_j>p_{j+1}$, and
the {\it descent number} $d_P$ of $P$ is the number of descents in it. Furthermore, given a word $P$ we shall denote by $C_P$ the
number \eqn\CP{ C_P:= \frac{(-1)^{d_P}}{|P|{|P|-1\choose d_P}}.  } For instance for $X = 25316$ we have $d_X=2$ and $C_X=
\frac{1}{5{4\choose 2} } =\frac{1}{30}$, while for $Y=351642$ we have $d_Y=3$ and $C_Y= -\frac{1}{6{5 \choose 3}}=
-\frac{1}{60}$.

The shuffle product between two
words is given by \Ree
\eqn\Shrecurs{
\emptyset\shuffle A = A\shuffle\emptyset = A,\qquad
A\shuffle B \equiv a_1(a_2 \ldots a_{n} \shuffle B) + b_1(b_2 \ldots b_{m}
\shuffle A)\,,
}
and it gives rise to all possible ways of interleaving the
letters of $A$ and $B$ without changing their original orderings within $A$ and $B$. For example
$12\shuffle 34 = 1234 + 1324 + 1342 + 3142 + 3124 + 3412$. The deconcatenation of a word $P$ into
two factors is denoted by $P=XY$ and it corresponds to all possible ways of splitting the word $P$ into
two words $X$ and $Y$. For example, if $P=123$ then $P=XY$ gives rise to the pair of words
$(X,Y)=(\emptyset,123),(1,23),(12,3),(123,\emptyset)$. The generalization to $P=X_1X_2\ldots X_k$ is straightforward.
Finally, the scalar product between two words $X$ and $Y$ is given by
\eqn\AdotB{
\langle X, Y\rangle \equiv \cases{$1$, & if $X=Y$;\cr
		 $0$, & otherwise.}
}
The deshuffle map $\d_k(P)$ is defined inductively as
\eqnn\deshuffledef
$$\eqalignno{
\d_k(i) &= \overbrace{i\otimes\emptyset\otimes\ldots\otimes\emptyset}^{k\rm\;times}
+ \emptyset\otimes i\otimes \ldots\otimes\emptyset
+ \cdots + \emptyset\otimes\emptyset\otimes \ldots\otimes i &\deshuffledef\cr
\d_k(i_1 \ldots i_n) &= \d_k(i_1)\cdots \d_k(i_n)\,,
}$$
or, equivalently, in terms of the shuffle product
\eqn\altdef{
\d_k(P) = \sum_{X_1, \ldots,X_k}\langle P,X_1\shuffle \ldots\shuffle X_k\rangle X_1\otimes \ldots\otimes X_k\,. 
}
The reduced deshuffle map $\bar{\d}_k(P)$ is obtained from $\d_k(P)$ by removing those terms which contain
the empty word $\emptyset$ as a tensor factor. 
For example, we have $\d_1(P) = P$,
$\d_2(12) = \emptyset\otimes 12  + 1\otimes 2 + 2\otimes 1 + 12\otimes\emptyset$, $\bar{\d}_2(12)=1\otimes 2+2\otimes 1$, $\bar{\d}_3(12)=0$ and
\eqn\exthree{
\d_2(123) =
  \emptyset\otimes123
+ 1\otimes23
+ 2\otimes13
+ 3\otimes12
+ 12\otimes3
+ 13\otimes2
+ 23\otimes1
+ 123\otimes\emptyset\,,
}
\eqn\exthreee{
\bar{\d}_2(123) =
  1\otimes23
+ 2\otimes13
+ 3\otimes12
+ 12\otimes3
+ 13\otimes2
+ 23\otimes1,
}
\eqn\exthreeee{
\bar{\d}_3(123) =
  1\otimes2\otimes3
+ 1\otimes 3\otimes2
+ 2\otimes1\otimes3
+ 2\otimes3\otimes1
+ 3\otimes1\otimes 2
+ 3\otimes2\otimes1.
}
We shall also adopt Sweedler's notation and write $\bar{\d}_k(P)=\sum_{(P)}P_{(1)}\otimes\ldots\otimes P_{(k)}$

We denote by $\min(P)$ the minimum among the letters in $P$.  We shall always deal with {\it multilinear words}, i.e., words
with no repeated letters. Given such a word $P$, its {\it standard factorization} $P=P_1\cdots P_k$ is defined iteratively as
follows. We put $j=\min(P)$ and consider the unique factorization $P=P'jP''$: if $P'=\emptyset$ we say that $P$ is a
(multilinear) {\it Lyndon word} and we define its standard factorization to be $P=P_1$, otherwise we take $P'=P_1\cdots P_{k-1}$
the standard factorization of $P'$, $P_k:=jP''$ and the standard factorization of $P$ is $P=P_1\cdots P_k$. Notice that by
construction all the factors $P_1,\ldots,P_k$ in the standard factorization of $P$ are Lyndon words, and $P_1>\cdots>P_k$ in the
lexicographical order: furthermore, the standard factorization is the only factorization of $P$ satisfying both these
properties\foot{In fact, this is a special case (the multilinear case) of the more general fact that any word admits a standard
factorization into a non-increasing product of Lyndon words, see \Reutenauer.}. For instance, the standard factorization of $X=
56427138$ is $X=X_1X_2X_3X_4 =(56)(4)(27)(138)$, and the standard factorization of $Y= 37528416$ is $Y=Y_1Y_2Y_3
=(375)(284)(16)$.

\newsec The color trace decomposition problem

Let us consider a simple Lie algebra whose generators $T^{a}$ satisfy\foot{In \Vcolor\ the
left-hand side is written in terms of a
representation-dependent normalization $I_{2R}$ as $\Tr(T_R^aT_R^b) = I_{2R}\d^{ab}$.
For convenience we use $I_{2R}=\half$ throughout this paper.}
\eqn\commut{
[T^a,T^b] = i f^{abc} T^c\,,\quad \Tr(T^a T^b) = \half\delta^{ab}\,,
}
where $f^{abc}$ denote the totally anti-symmetric structure constants and
$\d^{ab}$ is the Kronecker delta.
The symmetrized trace of Lie algebra generators
is defined by
\eqn\symtrace{
 d^{12\ldots n}\equiv {\rm Str}(T^1T^2 \ldots T^n) = {1\over n!} \, \sum_{\s \in S_n} \, \Tr
\big( T^{\s(1)} T^{\s(2)} \cdots T^{\s(n)}\big)\,,
}
where the sum is over all $n!$ elements of $S_n$
and we abbreviate the customary index of the Lie algebra generator $a_j$ simply by $j$.
Due to the
cyclicity of the trace we have from \symtrace\ and \commut\ that $d^{12}=\half \d^{12}$.

As discussed in \Vcolor, one is interested in decomposing the trace $\Tr(T^1 \ldots
T^n)$ in terms of symmetrized traces and structure constants leading to an expansion of the form
\eqn\stats{
\Tr(T^1 T^2 \ldots T^n) = d^{12 \ldots n} + \sum (fd + ffd + \cdots +
\underbrace{ff \cdots f}_{n{-}2\rm\;times}),
}
which can always be done in a systematic manner.
Using the algorithm implemented in the
{\tt color} package of {\tt FORM} \FORM\ and rewriting the results in the
color basis\foot{The results given by the
{\tt color} package are not written in a basis of color factors.} to be described below gives
the formulas in \Vtracetwo, with similar expansions at higher multiplicities (see e.g. the appendix B of \oneloopbb). 
These have been written in the basis of color factors
chosen in \oneloopbb, in which the letter $1$ is always in the symmetrized trace factor
$d^{1 \ldots}$. The remaining factors of structure constants will either have
contracted indices such as
$d^{12a}f^{45b}f^{b3a}$ and $d^{1ab}f^{24a}f^{35b}$ or will encompass all labels from $1$ to $n$ when no
factor of $d^{1 \ldots}$ is present. In both these cases we rearrange the
labels in such a way that the minimum and maximum labels are at the
extremities\foot{This choice is inspired by the del Duca--Dixon--Maltoni (DDM) basis \ddm.}
(where we consider a contracted index to be maximum).
This leads to basis elements such as $f^{13a}f^{a2b}f^{b45}$ or
$d^{12a}f^{34b}f^{b5a}$ and 
can be achieved using the Jacobi identities \Vcolor
\eqn\genJacobi{
f^{a[ij}f^{k]ab} = 0\,,\qquad d^{a(i_1i_2 \ldots i_{n-1}}f^{i_n)ab}  = 0\,.
}
As explained in \Vcolor, the decomposition \stats\ can always be done using
the following argument: starting from
the trivial identity
\eqn\trivid{
\Tr(T^1 \ldots T^n) = \Tr(T^1 \ldots T^n) - {\rm STr}(T^1 \ldots T^n) + d^{12 \ldots n}
}
one uses the commutation relation
\commut\ to move the generators in each one of the $n!$ terms in
${\rm STr}(T^1 \ldots T^n)$
to be in the same order as they appear in
$\Tr(T^1 \ldots T^n)$. Doing this for all $n!$ terms in $-{\rm Str}(T^1\ldots  T^n)$
cancels
the term $\Tr(T^1\ldots T^n)$ in the right-hand side of \trivid\
while generating lower-order terms containing structure constants as a result of the
commutation relation \commut\ and leading to \stats.

Before discussing the general solution to decomposing traces of color factors
we briefly review the definition of the Solomon idempotent.

\newsubsec\eid The Solomon idempotent

The {\it Solomon idempotent} appeared for the
first time
in the work of Solomon \solomon, who also noted its connection with Eulerian numbers. Hence the name {\it Eulerian idempotent}
is also commonly attributed to it\foot{Its characterization as a Lie idempotent was made by Reutenauer in \PBWReutenauer\ (see
also \Loday) but this aspect will not play a role in these notes.}. The Solomon idempotent appears
in several different contexts in the mathematical
literature such as in representations of
the symmetric group \refs{\garsia,\PBWReutenauer}, in free Lie algebras \Reutenauer, in
Hochschild homology \refs{\gersten,\giaquinto} and more recently
it has been used in connection with the Magnus series expansion solution
to differential equations \casas.

In order to define the Solomon idempotent, first recall the definition of the {\it descent number} $d_\sigma$ of the permutation $\sigma$,
\eqn\descno{
d_\sigma \equiv | \{ 1 \le i \le n-1 \, | \, \sigma(i) > \sigma(i+1)\}|\,.
}
For example, the permutation $43512$ has two descents (at the first and third
positions) so $d_{43512}=|\{1,3\}|= 2$.
In addition, we define
left-to-right nested commutators recursively by
$\ell(i_1,i_2, \ldots, i_n) \equiv
[\ell(i_1,i_2 \ldots,i_{n-1}),i_n]$, where $[i,j]=ij-ji$. For example $\ell(1,2,3,4)= [[[1,2],3],4]$. It is well known that
Lie polynomials with $n$ letters can be written in terms of the $(n{-}1)!$ dimensional {\it Dynkin basis}
$\ell(1,\s(2),\s(3),\ldots, \s(n))$.

The Solomon idempotent in the Dynkin basis of Lie polynomials is given by (this is shown in \bandiera, see also \casas)
\eqn\eulerianDynkin{
E(x_1x_2\cdots x_n) = \frac{1}{n}\sum_{\sigma\in S_n\atop\sigma(1)=1}
\frac{(-1)^{d_\sigma}}{{n-1 \choose d_\sigma}}\ell(x_1,x_{\sigma(2)}, \ldots ,x_{\sigma(n)})\,,
}
where $x_1, \ldots,x_n$ are non-commutative indeterminates. Using the notation from \CP, we might also write
\eqn\eulerianDynkinn{
E(x_1x_2\cdots x_n) =
\sum_{\sigma\in S_n\atop\sigma(1)=1}
C_{\sigma}\,\ell(x_1,x_{\sigma(2)}, \ldots ,x_{\sigma(n)})\,.
}
For instance for $n\le4$, defining $E^{12 \ldots n}\equiv E(T^1T^2 \cdots T^n)$ for the Solomon idempotent
written with Lie-algebra generators $T^j$,
formula \eulerianDynkin\ yields
\eqnn\fromeuler
$$\eqalignno{
E^1 &= T^1&\fromeuler\cr
E^{12} &= \half[T^1,T^2] \cr
E^{123} &=  {1 \over 3}\, {[ [ T^1 , T^2 ] , T^3 ]}
 -  {1 \over 6}\, {[ [ T^1 , T^3 ] , T^2 ]}\cr
E^{1234} &=
  {1 \over 4}\, {[ [ [ T^1 , T^2 ] , T^3 ] , T^4 ]}
 -  {1 \over 12}\, {[ [ [ T^1 , T^2 ] , T^4 ] , T^3 ]}
 -  {1 \over 12}\, {[ [ [ T^1 , T^3 ] , T^2 ] , T^4 ]}\cr
& -  {1 \over 12}\, {[ [ [ T^1 , T^3 ] , T^4 ] , T^2 ]}
 -  {1 \over 12}\, {[ [ [ T^1 , T^4 ] , T^2 ] , T^3 ]}
 +  {1 \over 12}\, {[ [ [ T^1 , T^4 ] , T^3 ] , T^2 ]}\cr
}$$
In view of \commut\ we define
\eqn\EPadef{
E^P := E^P_a T^a\,,
}
and note that, with the exception of
$E^1_a = \d^1_a$, all expansion coefficients $E^P_a$ are polynomials in the structure constants.
From $E^{12} = \half [T^1,T^2] = {i\over2}f^{12a}T^a$ we get
\eqn\Eextwo{
E^{12}_a={i\over2}f^{12a}.
}
And similarly,
\eqnn\Eexs
$$\eqalignno{
E^{123}_a &= -{1\over3}f^{12j}f^{j3a}+ {1\over6}f^{13j}f^{j2a},&\Eexs\cr
E^{1234}_a & =
-{i \over 4}\, f^{12j}f^{j3k}f^{k4a}
 +  {i \over 12}\,f^{12j}f^{j4k}f^{k3a}
 +  {i \over 12}\,f^{13j}f^{j2k}f^{k4a}\cr
&\quad{} +  {i \over 12}\,f^{13j}f^{j4k}f^{k2a}
 +  {i \over 12}\,f^{14j}f^{j2k}f^{k3a}
 -  {i \over 12}\,f^{14j}f^{j3k}f^{k2a}\,.
}$$
In general, using the notations \capF\ and \CP\ from the introduction, we may rewrite \eulerianDynkinn\ as
\eqn\eulerianexpansion{ 
E^{1\ldots n}_a =\sum_{\sigma\in S_n\atop\sigma(1)=1} i^{n-1}C_{\sigma} F^{\sigma}_a. 
}
A brief inspection of the expansions in \Vtracetwo\ and \Eexs\ reveals that the Solomon idempotent captures
the coefficients of the various terms in \Vtracetwo. This will be demonstrated below for the general case.

\newsubsec\PBWsolomon Trace decomposition from Solomon's projection

In order to obtain a closed formula that solves the color trace decomposition problem we
recall the projection obtained by Solomon in \solomon
\eqn\projS{
T^{P} = \sum_{k\ge1}\sum_{X_1,X_2, \ldots,X_k}{1\over k!}\langle
P, X_1\shuffle X_2\shuffle \ldots \shuffle X_k\rangle E^{X_1}E^{X_2} \ldots E^{X_k}
}
where $T^P \equiv T^{p_1}T^{p_2}\ldots T^{p_{n}}$ for a word $P=p_1p_2\ldots p_n$.
The multiplicity-two instance of \projS\ corresponds to the well-known decomposition into
a symmetric and antisymmetric combination (recall that $E^i= T^i$)
\eqn\twoex{
T^1T^2 = E^{12} + \half\big(E^{1}E^{2} + E^{2}E^{1}\big) = \half [T^1,T^2] + \half\big(T^{1}T^{2} + T^{2}T^{1}\big)\,.
}
But already at multiplicity three
\eqnn\solform
$$\eqalignno{
T^1T^2T^3 &= E^{123}
+ \half\big(
E^{12}T^{3}+E^{13}T^{2}+E^{23}T^{1}
+ T^{1}E^{23}+T^{2}E^{13}+T^{3}E^{12}
\big)&\solform\cr
&+{1\over3!}\big(T^{1}T^{2}T^{3} + T^{1}T^{3}T^{2} + T^{2}T^{1}T^{3} + T^{2}T^{3}T^{1} + T^{3}T^{1}T^{2}+ T^{3}T^{2}T^{1} \big)\,.
}$$
it is far from obvious that plugging in the expansions of the Solomon idempotents from \fromeuler\ into the right-hand
side recovers the monomial $T^1T^2T^3$ in the left-hand side.

As one can see from the above examples, the formula \projS\ projects
the product $T^1T^2 \ldots T^n$ into its totally symmetric component
${1\over n!}T^{(1}T^2 \ldots T^{n)}:=\frac{1}{n!}\sum_{\sigma\in S_n} T^{\sigma(1)}\ldots T^{\sigma(n)}$ plus lower-order terms containing
Eulerian idempotents. After taking the trace on both sides of Solomon's projection \projS,
the totally symmetric component is mapped to the symmetrized trace while the lower
order terms are mapped to sums of symmetrized traces multiplied by linear combinations
of structure constants as dictated by the Eulerian idempotents. This is the
solution to the color trace decomposition problem.

To see this more explicitly, we use the definition \EPadef\ to rewrite \projS\ as
\eqnn\projSA
$$\eqalignno{
T^P &=
\sum_{k\ge1}\sum_{X_1,X_2, \ldots,X_k}{1\over k!}\langle
P, X_1\shuffle X_2\shuffle \ldots \shuffle X_k\rangle E^{X_1}_{a_1}E^{X_2}_{a_2} \ldots
E^{X_k}_{a_k} T^{a_1}T^{a_2} \ldots T^{a_k}&\projSA\cr
& = \sum_{k\ge1}\sum_{X_1> X_2> \cdots>X_k}
\langle P, X_1\shuffle X_2\shuffle \ldots \shuffle X_k\rangle E^{X_1}_{a_1}E^{X_2}_{a_2} \ldots
E^{X_k}_{a_k} \tau^{a_1 \ldots a_k}
\,,
}$$
where we used that the shuffle product is
commutative to obtain the symmetrized product of the algebra generators by
ordering the sum according to $X_1>X_2> \cdots > X_k$ and defined $\tau^{a_1 \ldots a_k}={1\over k!}T^{(a_1} \ldots
T^{a_k)}$. Therefore multiplying \projSA\ by $T^1$ from the left, taking the trace
on both sides and using that $\Tr(T^{1}\tau^{a_1 \ldots a_k}) = d^{1 a_1\ldots a_k}$ leads to
\eqn\projSB{
\Tr{(T^{1P})} = \sum_{k\ge1}\sum_{X_1>X_2> \ldots > X_k}
\langle P, X_1\shuffle X_2\shuffle \ldots \shuffle X_k\rangle
E^{X_1}_{a_1}E^{X_2}_{a_2} \ldots E^{X_k}_{a_k} d^{1a_1a_2 \ldots a_k}
}
Alternatively, lifting the ordering restriction in the sum while compensating the overcount with ${1\over k!}$
and using \altdef\ leads to formula \closedshuffle\ from the introduction
\eqn\getclosedshuffle{
\Tr{(T^{1P})} = \sum_{k\ge1}\sum_{(P)}
{1\over k!} d^{1a_1a_2 \ldots a_k}\,
E^{P_{(1)}}_{a_1}E^{P_{(2)}}_{a_2} \ldots E^{P_{(k)}}_{a_k}\,,
}
concluding its proof.

Applying \getclosedshuffle\ for traces with up to six generators yields:
\eqnn\elegant
$$\eqalignno{
\Tr(T^1 T^2) &= d^{12} &\elegant\cr
\Tr(T^1 T^2 T^3) &= d^{123} +  d^{1a} E^{23}_a\cr
\Tr(T^1 T^2 T^3 T^4) &=d^{1234}   +  d^{12a} E^{34}_a + d^{13a} E^{24}_a
+d^{14a} E^{23}_a +  d^{1a}E^{234}_a\cr
\Tr(T^1 T^2 \ldots T^5) &=d^{12345} + d^{123a} E^{45}_a+d^{124a}
E^{35}_a+d^{125a} E^{34}_a+d^{134a} E^{25}_a+d^{135a} E^{24}_a \cr
&+d^{145a} E^{23}_a + d^{12a} E^{345}_a+ d^{13a} E^{245}_a+d^{14a} E^{235}_a+d^{15a} E^{234}_a \cr
&+ d^{1ab} (E^{23}_a E^{45}_b+ E^{24}_a E^{35}_b+ E^{25}_a E^{34}_b) + d^{1a}
E^{2345}_a \cr
\Tr(T^1 T^2 \ldots T^6) &=d^{123456} + d^{1234a} E^{56}_a+ d^{1235a} E^{46}_a + d^{1236a} E^{45}_a + d^{1245a} E^{36}_a + d^{1246a} E^{35}_a \cr
&+ d^{1256a} E^{34}_a + d^{1345a} E^{26}_a + d^{1346a} E^{25}_a + d^{1356a} E^{24}_a + d^{1456a} E^{23}_a \cr
&+ d^{123a} E^{456}_a+d^{124a} E^{356}_a+d^{125a} E^{346}_a+d^{126a} E^{345}_a+d^{134a} E^{256}_a \cr
&+ d^{135a} E^{246}_a+d^{136a} E^{245}_a+d^{145a} E^{236}_a+d^{146a} E^{235}_a+d^{156a} E^{234}_a \cr
&+ d^{12ab} (E^{34}_aE^{56}_b+E^{35}_aE^{46}_b+E^{36}_aE^{45}_b) \cr
&+ d^{13ab} (E^{24}_aE^{56}_b+E^{25}_aE^{46}_b+E^{26}_aE^{45}_b) \cr
&+ d^{14ab} (E^{23}_aE^{56}_b+E^{25}_aE^{36}_b+E^{26}_aE^{35}_b) \cr
&+ d^{15ab} (E^{23}_aE^{46}_b+E^{24}_aE^{36}_b+E^{26}_aE^{34}_b) \cr
&+ d^{16ab} (E^{23}_aE^{45}_b+E^{24}_aE^{35}_b+E^{25}_aE^{34}_b) \cr
&+ d^{12a} E^{3456}_a +d^{13a} E^{2456}_a +d^{14a} E^{2356}_a +d^{15a} E^{2346}_a +d^{16a} E^{2345}_a \cr
&+ d^{1ab} (E^{23}_aE^{456}_b+E^{24}_aE^{356}_b+E^{25}_aE^{346}_b+E^{26}_aE^{345}_b+E^{34}_aE^{256}_b \cr
&\quad{} +E^{35}_aE^{246}_b+E^{36}_aE^{245}_b+E^{45}_aE^{236}_b+E^{46}_aE^{235}_b+E^{56}_aE^{234}_b) \cr
& +  d^{1a} E^{23456}_a
}$$
It is not difficult to see that the
total number of terms generated by the formula \getclosedshuffle\ for
$n=3,4,5,6,7,8,9 \ldots$ is
equal to the Bell numbers $2,5,15,52,203,877,4140, \ldots$, respectively.

Finally, in order to obtain formula \eqclosedformula\ we look again at equation
\projSB, with the letter $1$ replaced by $0$ and the word $P$ replaced by $1\cdots n$.
For a fixed deshuffle $X_1\otimes\cdots\otimes X_k$ of $P$, we expand the
coefficients $E^{X_1}_{a_1}, \ldots, E^{X_k}_{a_k}$ in terms of structure constants
according to \eulerianDynkinn. Notice that for $1\leq j\leq k$ we have $X_j = i_1\ldots i_{|X_j|}$ with $i_1<\cdots< i_{|X_j|}$: by 
\eulerianexpansion\ we get $E^{X_j}_{a_j}=\sum_{\sigma_j} i^{|X_j|-1} C_{\sigma_j}F^{\sigma_j}_{a_j}$, where the sum runs over the words 
$\sigma_j$ obtained by permuting the last $|X_j|-1$ letters of $X_j$ while keeping the first one fixed. We obtain
\eqnn\eqend
$$\displaylines{
\Tr{(T^{01\ldots n})} = \hfil\eqend\hfilneg\cr
= \sum_{k\ge1}\sum_{X_1> \ldots > X_k}
\langle 1\ldots n, X_1\shuffle \ldots \shuffle X_k\rangle
\sum_{\sigma_1,\ldots,\sigma_k}i^{n-k}C_{\sigma_1}\ldots C_{\sigma_k} F^{\sigma_1}_{a_1}\ldots F^{\sigma_k}_{a_k} d^{0a_1a_2 \ldots a_k},
}$$
where the last sum runs over the words $\sigma_1,\ldots, \sigma_k$ obtained from $X_1,\ldots,X_k$ as above. Notice that each of
the words $\sigma_1,\ldots, \sigma_k$ has $\min(\sigma_j)=\min(X_j)$ as its first letter (in the terminology of subsection
\notsec, all the $\sigma_1,\ldots,\sigma_k$ are Lyndon words), and $\sigma_1>\cdots>\sigma_k$ in the lexicographical order
(since $X_1>\cdots>X_k$). Therefore, taking the permutation $\sigma\in S_n$ associated with the word $\sigma_1\cdots\sigma_k$,
the standard factorization of $\sigma$ as a word is precisely $\sigma=\sigma_1\cdots\sigma_k$. In the other direction,
given $\sigma\in S_n$ with standard factorization $\sigma=\sigma_1\cdots\sigma_k$,
we define $X_1,\ldots,X_k$ by rewriting the letters of $\sigma_1,\ldots,\sigma_k$
in increasing order: then $X_1\otimes\ldots\otimes X_k$ is a deshuffle of $1\ldots n$ with $X_1>\cdots> X_k$, and the words
$\sigma_1,\ldots,\sigma_k$ are obtained from $X_1,\ldots X_k$ as required in the last summation of \eqend.
This establishes a bijective correspondence between the terms of \eqend\ and \eqclosedformula, showing that the two formulas
are equivalent and concluding the proof of the latter.

\bigskip
\noindent{\bf Acknowledgements:}
We thank Oliver Schlotterer for collaboration on initial stages of this work and for sharing his notes
about patterns in the color trace decompositions.
CRM is supported by a University Research
Fellowship from the Royal Society.

\listrefs
\bye